\documentclass[10pt]{amsart}

\usepackage[square,compress,comma, numbers,sort]{natbib}
\usepackage[colorlinks=true, citecolor=blue, linkcolor=blue]{hyperref}
\usepackage{amsfonts}

\usepackage{amsfonts}
\allowdisplaybreaks[4]
\usepackage{color}

\def\IF{\infty}
\def\X{\vk X}
\definecolor{c22}{rgb}{0.2,0.6,0.3}
\definecolor{c20}{rgb}{0.,0.7,0.}
\definecolor{c30}{rgb}{0.,0.,1.}
\definecolor{c40}{rgb}{1,0.1,0.7}
\definecolor{c50}{rgb}{1,0,0}
\definecolor{c60}{rgb}{1,0.9,0.1}

\def\zE#1{\textcolor{c20}{#1}}
\def\tE#1{\textcolor{c20}{#1}}
\def\sE#1{\textcolor{c22}{#1}}
\def\sE#1{#1}
\def\pE#1{\textcolor{c20}{#1}}
\def\pE#1{#1}
\def\zE#1{#1}
\def\tE#1{#1}

\def\cE#1{\textcolor{c20}{#1}}

\def\cE#1{#1}

\def\cJ#1{\textcolor{c50}{#1}}
\def\cJ#1{#1}

\def\mE#1{\textcolor{c20}{#1}}
\def\mE#1{#1}

\def\kE#1{{\textcolor{c20}{#1}}}
\def\kE#1{#1}

\def\cL#1{\textcolor{c30}{#1}}
\def\cL#1{#1}
\def\cLL#1{\textcolor{c40}{#1}}
\def\cLL#1{#1}
\def\ccL#1{\textcolor{c50}{#1}}
\def\ccL#1{#1}
\def\KK#1{{\textcolor{c20}{#1}}}
\def\KK#1{#1}

\def\cJi#1{\textcolor{c50}{#1}}
\def\cJi#1{#1}

\def\ccLP#1{\textcolor{c50}{#1}}
\def\ccLP#1{#1}

\def\cccL#1{\textcolor{c50}{#1}}
\def\cccL#1{#1}
\def\ccJ#1{\textcolor{c50}{#1}}
\def\ccJ#1{#1}

\newcommand{\kb}[1]{\boldsymbol{#1}}
\newcommand{\vk}[1]{\kb{#1}}

\newcommand{\ve}{\varepsilon}
\def\fracl#1#2{\left(\frac{#1}{#2} \right) }

\newcommand{\abs}[1]{\left\lvert #1 \right\rvert}
\newcommand{\Abs}[1]{ \biggl \lvert #1 \biggr \rvert}

\newcommand{\E}[1]{\mathbb{E}\left\{#1\right\}}

\newcommand{\pk}[1]{\mathbb{P} \left\{ #1 \right\} }

\newcommand{\R}{\mathbb{R}}

\newcommand{\N}{\mathbb{N}}
\newcommand{\inr}{\in \R}
\newcommand{\inn}{\in \N}
\newcommand{\ldot}{,\ldots,}

\newcommand{\limit}[1]{\lim_{#1 \to   \infty}}

\newcommand{\BQN}{\begin{eqnarray}}
\newcommand{\EQN}{\end{eqnarray}}
\newcommand{\BQNY}{\begin{eqnarray*}}
\newcommand{\EQNY}{\end{eqnarray*}}

\newcommand{\BS}{\begin{sat}}
\newcommand{\ES}{\end{sat}}
\newcommand{\BT}{\begin{theo}}
\newcommand{\ET}{\end{theo}}
\newcommand{\BK}{\begin{korr}}
\newcommand{\EK}{\end{korr}}

\newcommand{\BD}{\begin{de}}
\newcommand{\ED}{\end{de}}

\newcommand{\BIT}{\begin{itemize}}
\newcommand{\EIT}{\end{itemize}}
\newcommand{\BDI}{\begin{description}}
\newcommand{\EDI}{\end{description}}

\newcommand{\BRM}{\begin{remarks}}
\newcommand{\ERM}{\end{remarks}}

\newcommand{\BEL}{\begin{lem}}
\newcommand{\EEL}{\end{lem}}

\newtheorem{theo}{Theorem}[section]
\newtheorem{sat}[theo]{Proposition}
\newtheorem{de}[theo]{Definition}
\newtheorem{lem}[theo]{Lemma}

\newtheorem{korr}[theo]{Corollary}

\newtheorem{remarks}[theo]{Remarks}

\newcommand{\nelem}[1]{{Lemma \ref{#1}}}

\newcommand{\netheo}[1]{{Theorem \ref{#1}}}

\newcommand{\nesec}[1]{{Section \ref{#1}}}

\newcommand{\prooftheo}[1]{ \textsc{\bf Proof of Theorem} \ref{#1}:}

\newcommand{\prooflem}[1]{\textsc{\bf Proof of Lemma} \ref{#1}:}

\newcommand{\COM}[1]{}

\newcommand{\QED}{\hfill $\Box$}

\def\cE#1{#1}

\def\IF{\infty}

\newcommand{\expon}[1]{\exp\left(#1\right)}

\def\I#1{\mathbb{I}\{#1\}}

 \def\Y{\pE{\vk{Y}}}

\def\Cov{\mathrm{Cov}}
\newcommand{\toprob}{ \stackrel{p}{\to}}
\newcommand{\equaldis}{\stackrel{d}{=}}

\newcommand{\todis}{\stackrel{d}{\to}}

\begin{document}
\title{Extremes and Limit Theorems  for Difference of Chi-type processes}\thanks{E. Hashorva, L. Ji and C. Ling were  partially  supported by the Swiss National Science Foundation (SNSF) grant 200021-140633/1, E. Hashorva is partially supported from SNSF grant 200021-166274, C. Ling also acknowledges Fundamental Research Funds for Central Universities (XDJK2016C118), and Fundamental Science and Advanced Technology Funds of Chongqing, China (cstc2016jcyjA0036).}
\author{Patrik Albin}\address{Patrik Albin, Department of Mathematical Sciences, Chalmers University of
Technology, \cccL{SE-412, 96 Gothenburg}, Sweden}
\author{Enkelejd  Hashorva}\address{Enkelejd  Hashorva, Department of Actuarial Science,
University of Lausanne, \cccL{UNIL-Dorigny, 1015  Lausanne,} Switzerland}
\author{Lanpeng Ji}
\address{Lanpeng Ji, Department of Actuarial Science,
	University of Lausanne, \cccL{UNIL-Dorigny, 1015  Lausanne,} Switzerland}
\author{Chengxiu Ling}
\address{School of Mathematics and Statistics, Southwest University, \cccL{Beibei District, 400715 Chongqing,} China}
%
%
\begin{abstract}
Let $\{\zeta_{m,k}^{(\kappa)}(t), t \ge0\},
\kappa>0$ be {random processes} defined as the
differences of  two independent stationary chi-type processes with
$m$ and $k$ degrees of freedom.
%
 In this paper we derive the asymptotics of $\pk{\sup_{t\in[0,
T]}\zeta_{m,k}^{(\kappa)}(t)> u}, {u\to\IF}$  under some
assumptions on the covariance structures of the underlying Gaussian
processes. Further, we establish a Berman sojourn limit theorem  and a Gumbel limit \kE{result}.
  \end{abstract}

\maketitle
\section{Introduction} \label{sec1}
\kE{Let}  $\vk{X}(t)=(X_1(t),\ldots,X_{m+k}(t)), t\ge0, m\ge 1, k\ge 0$ be a   vector process with independent components which are
centered stationary Gaussian processes with almost surely {(a.s.)}
continuous sample paths, \mE{and unit variance. Set  $r_i(t)=\E{X_i(t) X_i(0)},t\ge 0$ and suppose that }
\BQN\label{corrr}
 r_i(t)= 1 -  \cE{C_i}\abs{t}^{\alpha}+ o(\abs{t}^{\alpha}), \quad  t\to 0\quad {\rm and}\quad r_i(t)< 1, \quad  \forall t{\ne}0,\ \cL{1\le i\le m+k,}
\EQN
where $\alpha \in (0,2]$ \cE{and {$\vk C:=(C_1,
\ldots, C_{m+k})\in \tE{(0,\IF)^{m+k}}$}.} Define \mE{the random process\cccL{es}} $\left\{\zeta_{m,k}^{(\kappa)}(t),   t\ge0
\right\}$, $\kappa>0$ by
\BQN \label{ModelA}
\zeta_{m,k}^{(\kappa)}(t):= \Biggl( \sum_{i=1}^m X_i^2(t)
\Biggr)^{\kappa/2}- \Biggl(\sum_{i=m+1}^{m+k}
\cE{X_i^2(t)}\Biggr)^{\kappa/2}=:|\X^{(1)}(t)|^{\kappa} -
 |\X^{(2)}(t)|^{\kappa} , \quad t\ge0.
\EQN
In this paper we shall investigate \kE{the asymptotics of}
\BQN
\label{eq:zeta} \pk{\sup_{t\in[0, T]} \zeta_{m,k}^{(\kappa)}(t) >u},
\quad u\to\IF
\EQN
\mE{for given positive  $T$}. 
 \mE{The study of the asymptotics of \eqref{eq:zeta} is of interest  in
  engineering applications concerned} with the  safety of structures;  
  see, e.g.,
\cite{Lindgren1980,Lindgren1984,Lindgren1989} and the references
therein.  \ccJ{More \mE{precisely}},   of interest is the \mE{calculation of the} probability that  \cJi{the Gaussian vector
process} {exits a predefined \KK{safety}  region}
$\vk{S}_u\subset\R^{\sE{m+k}}$ up to the time $T$, \mE{i.e., }
\BQNY
\pk{\vk{X}(t)\not\in\vk{S}_u,\ \ \text{for}\ \text{some}\
t\in[0,T]}. \EQNY
\mE{In the
aforementioned papers}, various types of \KK{safety}  regions $\vk{S}_u$
\ccJ{have been} considered for smooth Gaussian \sE{vector} processes. Particularly, a  {safety} region given by a
ball centered at 0 with radius $u>0$
$$
\vk{B}_u =\Biggl\{(x_1,\ldots,x_{m+k})\in\R^{\sE{m+k}}:
\Biggl(\sum_{i=1}^{m+k} x_i^2\Biggr)^{1/2}\le u\Biggr\}
$$
has been {extensively studied}; see,  \kE{e.g.,}  \cite{Albin1990, Berman82,
Claudia2013, Pit96}. Referring to \cite{AlbinPHD, Albin1990},  we
know that  for $k=0$ 
\BQNY
\pk{\vk{X}(t)\not\in\vk{B}_u,\ \ \text{for}\ \text{some}\
t\in[0,T]}&=& \pk{\sup_{t\in[0, T]} \abs{\vk X(t)} >u}
\nonumber\\
&=& T H_{\alpha,1}^{m,0}(\vk C) u^{\frac{2}{\alpha}}
\pk{\abs{\vk X(0)} >u}(1+o(1)), \quad u\to\IF,
\EQNY
where $\tE{H_{\alpha,1}^{m,0}(\vk C)}$ is a positive constant
(see \eqref{Pickands constant} below for a precise definition).  Very recently \cite{Zhdanov} obtained the tail asymptotics of
 the product of two Gaussian processes \ccJ{which has the same tail asymptotic behavior as   $ \sup_{t\in [0,T]} \zeta_{1,1}^{(2)}(t)$.}\\
 Our first result extends the findings of \kE{\cite{Albin1990,Zhdanov}} and suggests   an asymptotic
 approximation for the exit probability of $\vk X$ from the {safety} regions $ \vk{S}_u^{(\kappa)}$  defined by
\BQNY\label{def.Su} \vk{S}_u^{(\kappa)} = \left\{ (x_1, \ldots,
x_{m+k})\in \R^{m+k}: |\vk x^{(1)}|^{\kappa}-
|\vk x^{(2)}|^{\kappa}\le u\right \}.
\EQNY
 Since  chi-type processes  appear  naturally
as limiting processes  \ccJ{(see\cL, e.g., \cite{Albin2003, Aue2009, LingT2015})},  when
{one considers} two independent asymptotic models, the study of the
supremum of the difference of the two chi-type processes  \tE{is} \mE{also} of
some interest in mathematical statistics
and its applications.

Although for $k\ge1$ the random process  $\zeta_{m,k}^{(\kappa)} $
is not Gaussian   and the analysis of \KK{the} supremum \mE{cannot be directly}
{transformed} into the study of {the} supremum of a related Gaussian
random field (which is the case for  chi-type processes; see, e.g.,
\cite{HashPengChiTrend,Lindgren1980,Lindgren1984,Lindgren1989, LingP2016, 
Pit96, TanHash13}), it turns out that it is possible to apply the
techniques for dealing with extremes of stationary
processes developed mainly in \cite{Albin1990,Berman82,Berman92}.  In the second part of Section 2 we
derive a sojourn limit theorem for $\zeta_{m,k}^{(\kappa)}$.
\cJ{Further, we \KK{show} a Gumbel limit theorem for the supremum of
$\zeta_{m,k}^{(\kappa)}$ over an increasing infinite  interval. We
refer to \cite{Albin1990, Albin2003, Aue2009, leadbetter1983extremes, Pit96,tanH2012} for results on the Gumbel
limit theorem for \KK{Gaussian} processes and chi-type processes.}

Brief outline of the paper: \cL{o}ur main results are stated in Section
2. In \nesec{sec3} we present  proofs of \netheo{T1}, \netheo{T2}
and \netheo{T3} followed then by an appendix.

\section{Main Results}
We first introduce some  notation. Let  $\{Z(t), t \ge 0\}$ be
a standard fractional Brownian motion (fBm) with Hurst index
$\alpha/2 \in (0,1]$, i.e., it is a centered Gaussian process with
a.s. continuous sample paths and covariance function
$$
 \Cov(Z(s), Z(t))= \frac12 \Bigl( s^\alpha+ t^\alpha- \abs{s-t}^\alpha\Bigr), \quad s,\, t \ge 0.
$$
In the following, \cJi{let} $\{Z_i(t), t \ge 0\}, 1\le i\le m+k$ be
independent
copies of $Z$ and  define  $\mathcal{W}_\kappa$ to be a  Gamma
distributed random variable with parameter  $ {(k/\kappa, 1)}$. Further let
$\vk{O}_1 = (O_{1}, \ldots, O_{m}), \vk{O}_2= (O_{m+1}, \ldots,
O_{m+k})$ denote two random vectors uniformly distributed on the
unit sphere of $\R^{{m}}$ and $\R^k$, respectively. Hereafter we
shall suppose that $ \vk{O}_1, \vk{O}_2,\mathcal{W}_\kappa $ and  $Z_i$'s
are mutually independent. Define for $ m\ge1, \ k\ge0,\  \ccL{\kappa>0}$
\BQN\label{eq:eta} \eta_{m,k}^{(\kappa)}(t)= {\widetilde Z_{m,k}^{\ccLP{(\kappa)}}(t)
+E,}\quad t\ge0, \EQN where $E$ is a unit mean exponential random
variable being independent of all {the} other random elements
involved, and \ccJ{(recall $\vk C=(C_1,\ldots, C_{m+k})$ given in \eqref{corrr})} 
\BQN \label {KL}
\widetilde Z_{m,k}^{\ccLP{(\kappa)}}(t) = \left\{
\begin{array}{ll}
L_1(t),&\kappa>1,\\
L_1(t)+L_2(t),&\kappa=1,\\
L_2(t),&\kappa<1,
\end{array}
\right.
\mbox{with}\ \left\{
\begin{array}{l}
L_1(t)=\sum_{i =
1}^m \sqrt{2C_i} O_iZ_i(t) -\left( \sum_{i=1}^m C_i
O^2_{i}\right) t^\alpha,\\
L_2(t)=\mathcal W_\kappa -\left(\mathcal
W_\kappa^{2/\kappa}+2{(\mathcal W_\kappa/\kappa)}^{1/\kappa}\sum_{i=m+1}^{m+k} \sqrt{2C_i}O_iZ_i(t)
\right. \\
\left.\ \ \qquad\quad +2\kappa^{-2/\kappa}\sum_{i=m+1}^{m+k}C_iZ_i^2(t)  \right)^{\kappa/2},
\end{array}
\right.
\EQN
Here we set $\sum_{i=m+1}^mc_i =:0$.\\
We state next our main  result.
\BT \label{T1}
\cJi{If $\{\zeta_{m,k}^{(\kappa)}(t), t\ge0\}$
is given by
\eqref{ModelA} with the involved Gaussian processes   $X_i$'s  satisfying} \eqref{corrr},  \ccLP{then, for any  $T>0$ }
\BQNY\label{eq:main}
\pk{ \sup_{t \in [0,T]} \zeta_{m,k}^{(\kappa)}(t)> u}
 = T H^{m,k}_{\alpha,\kappa} (\vk C)  u^{\cccL{\frac{2\tau}{\alpha\kappa}}}  
 \pk{\zeta_{m,k}^{(\kappa)}(0)> u}  (1+o(1))
\EQNY
\ccJ{holds  as $u\to\IF$,}
where \cccL{$\tau=2/\kappa -1$ for $\kappa\in(0,1)$, and 1 otherwise,} 
 with \ccLP{$\eta_{m,k}^{(\kappa)}$ given by \eqref{eq:eta},}
\BQN\label{Pickands constant}
H^{m,k}_{\alpha,\kappa} (\vk C) = \lim_{a\downarrow 0}\frac 1 a
\pk{\sup_{j\ge 1}\eta_{m,k}^{(\kappa)}(aj)\le 0} \in (0,\IF).
\EQN
\ET
 {\bf Remarks:} a) \ccJ{The tail asymptotics of the Gaussian chaos $\zeta_{m,k}^{(\kappa)}(0)$ is discussed in Lemma \ref{LTA} below.}

 \ccJ{b) The most obvious choice of $\kappa$ is 1, which corresponds to the difference of $L_2$-norm of two independent multivariate Gaussian processes. For the case $\kappa=2$ and $m=k=1$ the problem was (implicitly) investigated by considering the product of two independent Gaussian processes in the recent contribution \cite{Zhdanov}.}

 c) \ccJ{Since $\vk O_1$ is uniformly distributed on the
unit sphere of $\R^{{m}}$, we have, for $\kappa>1$ and $\vk C=\vk 1$, that $\eta_{m,k}^{(\kappa)}(t) \stackrel d= \sqrt2Z(t)-t^\alpha +E$. In such a case, the constant $H^{m,k}_{\alpha,\kappa}(\mathbf 1)$ coincides with  the classical Pickands constant $H_\alpha$; see, e.g., \cite{Pit96}. Approximation of Pickands constant $H_\alpha$ has been considered by a number of authors; see
 the recent contribution \cite{DikerY} which gives some simulation algorithms. Precise estimation of the general  Pickands constant $H^{m,k}_{\alpha,\kappa}(\vk C)$  seems to be hard to find, due to the complexity of the process $\eta_{m,k}^{(\kappa)}$.}

 \COM{The remaining cases are more involved due to the polar representation of Gaussian vector $\vk X$ and the value $\kappa\le 1$, see \nelem{LTA2}. The complicated  processes $\eta$ appear intuitively according to the three tail behavior of $\zeta_{m,k}^{(\kappa)}$, see \nelem{LTA}. Note that  the recent contribution  {\cite{DikerY}
 {gives} simulation algorithms } that {throw}
 light on Pickands constant  $H_\alpha$.}

\COM{
c) Define exit times $\tau_{\kappa}(u)=\inf\{t>0: \vk X(t)\not\in
\vk S_u^{(\kappa)}\}, \kappa>0$ with $\vk S_u^{(\kappa)}$ given by
\eqref{def.Su}. \tE{By a direct application of \netheo{T1} for}
any $T>0$ we obtain
\BQNY \limit{u} \pk{\tau_{\kappa}(u)\le t|
\tau_{\kappa}(u)\le T} = \frac tT,\quad \forall t\in[0, T],
 \EQNY
\kE{which} means that asymptotically $\tau_{\kappa}(u)|\{\tau_{\kappa}(u)\le T\}$ is uniformly distributed on $[0, T]$.\\
}
d) 
 We see from \netheo{T1} and Lemma \ref{LTA} that, if $\kappa>2$,   
 then, for any $m,k\ge 1$
\BQNY
\pk{ \sup_{t \in [0,T]} \zeta_{m,k}^{(\kappa)}(t)> u} = \pk{ \sup_{t \in [0,T]} \zeta_{m,0}^{(\kappa)}(t)> u} (1+o(1))
\EQNY
\ccJ{holds as $u\to \IF$, which means that
  $X_{m+1}\ldot X_{m+k}$ do not influence the tail asymptotic of $\sup_{t \in [0,T]} \zeta_{m,k}^{(\kappa)}\ccLP{(t)}$. This is not so surprising as the tail asymptotic behavior of $\zeta_{m,k}^{(\kappa)}(0)$ is subexponential.}

\ccJ{Next, we} consider the sojourn time of $ \zeta_{m,k}^{(\kappa)}$ above a threshold
$u
> 0$ in the time interval $[0,t]$ defined by
\BQNY
L_{m,k,t}^{(\kappa)}(u) = \int_0^t
\I{\zeta_{m,k}^{(\kappa)}(s) > u} \, ds, \quad t>0, 
 \EQNY
 \cccL{where $\I{\cdot}$ stands for the indicator function.} 
\kE{Our second result below establishes a Berman sojourn limit theorem for $\zeta_{m,k}^{(\kappa)}$. } {See \cite{Berman82} for related discussions on sojourn times of Gaussian processes and related processes.}
 \BT \label{T2}Under the assumptions and notation of \netheo{T1}, we have, for
any $t>0$
\BQNY
\int_x^\IF \pk{u^{\ccL{\frac{2\tau}{\alpha\kappa}}}L_{m,k,t}^{(\kappa)}(u)> y  }\, dy =
 u^{\cccL{\frac{2\tau}{\alpha\kappa}}} \E{ L_{m,k,t}^{(\kappa)}(u)}
\Upsilon_{\kappa}(x) (1+o(1)) 
 \EQNY
 holds \ccJ{as $u\to\IF$} for all continuity point  $x>0$ of  $\Upsilon_{\kappa}(x) := \pk{\int_0^\IF \I{\eta_{m,k}^{(\kappa)}(s)> 0}\, ds > x}$. 
\ET

\COM{  {\bf Remarks}: a) It {might be}  possible to allow $X_i$'s to be dependent. Results for
 extremes of chi-type processes for such generalizations can be found in \cite{Albin2003,Aue2009}.\\
b) Following the methodology in \cite{Albin94}
{one could consider} $X_i$'s to be self-similar Gaussian
processes.
Further extensions for {random} fields \ccLP{could also be}
possible by adopting the recent findings in \cite{MR2733939}. \\}

In the following, we  derive a Gumbel limit theorem  for
$\sup_{t\in [0, T]} \zeta_{m,k}^{(\kappa)}(t)$ under a linear
normalization, which is also of interest in extreme value analysis and statistic\cL{al} test\cL{s}. \ccJ{ We refer to \cite{Berman92,leadbetter1983extremes, Aue2009, Jaruvskova2015detecting}  for its applications in deriving approximations of  the critical values of  the proposed test statistics.}
\COM{For \zE{simplicity} we assume that all of the
covariance functions {$r_i(\cdot), 1\le i\le m+k$  are  the same,
denoted by $r(\cdot)$}. Hence we have that condition \eqref{corrr}
holds with $C_i=C>0, 1\le i\le m+k$.}

\BT \label{T3}
\kE{Under the assumptions and notation of \netheo{T1},
if further the following Berman\cL{-type} condition}
\BQN\label{Berman-type}
\limit{t} {\max_{1\le l\le m+k}|r_l(t)| }\cL{(\ln t )^c}=0, \ \ \text{with}\ c=\left\{\begin{array}{ll}
2/\kappa-1,&0<\kappa<1,\\
1,&1\le \kappa\le 2,\\
k+1-2k/\kappa,&\kappa>2
 \end{array}\right.
\EQN
holds, 
then
\BQNY
 \limit{T}\sup_{x\in \R} \Abs{ \pk{ a_T^{(\kappa)} \Big(\sup_{ t\in [0, T]}\zeta_{m,k}^{(\kappa)}(t)- b_T^{(\kappa)}\Big) \le x} -
\expon{ - e^{-x}}} =0,
\EQNY
where, for all $T$ large
\BQN \label{Normalized
constant}
 \quad a_T^{(\kappa)} = \frac{(2 \pE{\ln } T)^{1-\kappa/2}}\kappa,
 \quad b_T^{(\kappa)} = (2 \pE{\ln } T)^{\kappa/2} + \frac\kappa{2(2 \pE{\ln } T)^{1-\kappa/2}} \left(
K_0\ln\ln T+
 \pE{\ln } D_0\right),
 \EQN
 with \mE{(below $\Gamma(\cdot)$ denotes the Euler Gamma function)}
  \begin{align*}
 \cL{D_0}&=\left(\frac{H_{\alpha,\kappa}^{m,k}(\vk C)
}{\Gamma(m/2)\Gamma(k/2)}\right)^2\times \left\{
 \begin{array}{ll}
 2^{\frac{{2}}\alpha(\frac{2}{\kappa}-1) +2\left(1-\frac k\kappa\right)}
 \left(\Gamma\left(\frac k\kappa\right)\kappa^{\ccL{(k/\kappa-1)}}\right)^{2},&0<\kappa\le 1,\\
 2^{\frac{{2}}\alpha +2\left(1-\frac k\kappa\right)}
 \left(\Gamma\left(\frac k\kappa\right)\kappa^{\ccL{(k/\kappa-1)}}\right)^{2},&1<\kappa<2,\\
 2^{\frac{{2}}\alpha-\cLL{2}}\left(\Gamma\left(\frac k2\right)\right)^{2},&\kappa=2, \\
  2^{\frac{{2}}\alpha}\left(\Gamma\left(\frac k2\right) \right)^{2},&\kappa>2
  \end{array}
 \right.
\\
K_0&=\left\{
\begin{array}{ll}
m-2+(2/\alpha)(2/\kappa-1)+ k(1-2/\kappa), & 0<\kappa \le 1,\\
m-2+2/\alpha+ k(1-2/\kappa), & 1<\kappa <2,\\
m-2+2/\alpha, & \kappa \ge2.
\end{array}
\right.
\end{align*}
 \ET

 \COM{ {\bf Remar\cL{ks}:}  In contrast to the Normal Comparison Lemma in  \cite{leadbetter1983extremes} for the Gaussian processes, we obtain the crucial comparison inequality of the chi-type processes in \nelem{LTD1} to ensure the mixing D condition, which was
first addressed by the seminal paper \cite{LeadbetterR1982}  (see
Lemma 3.5 therein and also Theorem 10 in \cite{Albin1990}).
\nelem{LTD1} is also  expected to be useful in extreme value analysis {when} concerned with chi-type processes; see, e.g., \cite{tanH2012}  avoiding the technical verification of mixed-Gumbel limit theorems for the strongly dependent cyclo-stationary $\chi$-processes . \\
}

Under the assumptions of \netheo{T3},  we have the following convergence in probability (denoted by
$\toprob$)
 \BQNY \frac{\sup_{ t\in [0, T]}\zeta_{m,k}^{(\kappa)}(t)}{(2
\ln T)^{\kappa/2}} \toprob 1, \quad {T}\to \IF,
 \EQNY
 \cJi{which
\kE{follows} from the fact that $\lim_{T\to\IF}b^{(\kappa)}_T / (2 \ln
T)^{\kappa/2}= 1$ and  that $a_T^{(\kappa)}$ is bounded away from zero,
together with elementary considerations.}
 In several \cJ{cases} such a convergence in probability can be strengthened
 to the {$p$th} mean convergence which is referred to as the Seleznjev $p$th mean convergence
 since the idea \ccJ{was first suggested by Seleznjev}  in \cite{MR2275916}, see also \cite{MR2023892}.   In order to show  the Seleznjev $p$th mean convergence
 \ccL{of crucial importance} is the Piterbarg inequality (see \cite{Pit96}, Theorem 8.1).  Since the Piterbarg inequality holds also for chi-square processes  (see {\cite{TanHash13}, Proposition 3.2}),  using \ccL{further} the fact that
$$ \zeta_{m,k}^{(\kappa)}(t) \le |\X^{(1)}(t)|^{\kappa},
\quad t\ge 0,
$$
we immediately get the Piterbarg inequality for the difference  of
  chi-type processes by simply applying the aforementioned
proposition. Specifically, under the assumptions of  \netheo{T3} for
any $T>0$ and all large $u$
\BQNY \pk{\sup_{t \in
[0,T]}\zeta_{m,k}^{(\kappa)}(t)> u} \le K T u^{\beta}
\expon{-\frac{1}{2}u^{2/\kappa}},
\EQNY
 where $K$ and $\cJ{\beta}$
are two positive constants not depending on $T$ and $u$. Note that
the above {result} also follows immediately from \netheo{T1}
\cJi{combined with \nelem{LTA} below.} Hence utilizing Lemma 4.5 in
\cite{TanHash13} we arrive at our last result.

\BK (Seleznjev $p$th mean theorem) Under the assumptions of
\netheo{T3}, we have, for any $p>0$
\BQNY
\limit{T}\E{\fracl{\sup_{
t\in [0, T]}\zeta_{m,k}^{(\kappa)}(t)}{(2 \ln T)^{\kappa/2}}^p }= 1.
\EQNY
\EK
\section{Further Results and Proofs}
\label{sec3}
\ccJ{Before presenting the proof of \netheo{T1} we  first}  give some preliminary lemmas. Hereafter we use the
same notation and  assumptions as in \nesec{sec1}. By $\todis$ \ccL{and $\equaldis$} we
shall denote the convergence in distribution (or the convergence of
finite dimensional distributions if both sides of it are random
processes) \ccL{and equality in distribution function, respectively}.   Further, we write $f_\xi{(\cdot)}$ for the pdf of a
random variable $\xi$ and {write} $h_1\sim h_2$ if  two functions
$h_i(\cdot), i=1,2$ {are}  such that $h_1/h_2$ goes to 1 as the
argument tends to some limit.  {For simplicity} we shall denote, with
$\kappa>0$ and $\tau=\max(2/\kappa-1, 1)$,
$$q_\kappa = q_\kappa(u) = u^{-2{\tau}/(\alpha\kappa)}, \quad w_\kappa(u)=\frac1\kappa u^{2/\kappa-1}, \quad u>0.$$
In the proofs of Lemmas \ref{LTA}--\ref{LTB}, we denote
$u_{\kappa,x}=u+x/w_\kappa(u)$ for all $u, x>0$.

\BEL \label{LTA} {Let $\{\zeta_{m,k}^{(\kappa)}(t), \ t\ge0\}$ be given by \eqref{ModelA}}.  For all integers $m\ge1, k\ge0$ we have as $u\to\IF$
\begin{align*} \pk{ \zeta_{m,k}^{(\kappa)}(0) > u} \sim \frac{f_{\zeta_{m,k}^{(\kappa)}(0)}(u)}{w_\kappa(u)}
\sim
\frac{2^{2-(m+k)/2}}{\kappa^2\Gamma(k/2)\Gamma(m/2)} \frac{u^{m/\kappa-1}}{w_\kappa(u)} \expon{-\frac12 u^{2/\kappa}}
\cccL{\times}\left\{\begin{array}{ll}
\frac{\Gamma(k/\kappa)}{(w_\kappa(u))^{k/\kappa}}, &\kappa<2,\\
\Gamma(k/2), & \kappa=2,\\
\kappa 2^{k/2-1}\Gamma(k/2), & \kappa>2,
\end{array}
\right.
\end{align*}
where $\Gamma(k/ \kappa) / \Gamma(k/2):=1$ for $k=0$ and all $\kappa>0$.
 \EEL

{\bf Proof.} \mE{The claim follows from \cL{Theorem 1 in \cite{Korshunov}}}.
\QED

\BEL \label{LTA2}
If  $\{\zeta_{m,k}^{(\kappa)}(t),t\ge0\}$ is as in \netheo{T1}, then

\BQNY
\left\{ w_\kappa(u)(\zeta_{m,k}^{(\kappa)}(q_\kappa
t)- u)  \lvert \{\zeta_{m,k}^{(\kappa)}(0)>u\},\ t\ge0\right\} \todis
\left\{\eta_{m,k}^{(\kappa)}(t),\ t\ge0\right\}, \quad u\to\IF,
 \EQNY
  with
  $\eta_{m,k}^{(\kappa)}$ given by  \eqref{eq:eta}. {Recall that $\todis$ stands for the convergence of finite dimensional distributions}.
 \EEL
{\bf Proof.} \cJi{We henceforth adopt the notation introduced in Section 2.}
\kE{By \nelem{LTA}}\cL{, we have}
\BQNY
 w_\kappa(u)( \zeta_{m,k}^{(\kappa)}(0)- u) \Bigl\lvert \{\zeta_{m,k}^{(\kappa)}(0)> u\}\ {\todis}\   E, \quad u\to \IF.
\EQNY
 Thus, in view of  Theorem 5.1 in \cite{Berman82},  it suffices to show that, for any $0 < t_1 < \cdots < t_n  < \IF, \cJi{n\inn}$
\begin{align} \label{Proof_Z}
  p_{k}(u)&:=\pk{\cap_{j=1}^n\{\zeta_{m,k}^{(\kappa)}(q_\kappa t_j)\le u_{\kappa, z_j}\}\Bigl
\lvert \zeta_{m,k}^{(\kappa)}(0)= {u_{\kappa,x}}}\nonumber\\
 &\to
\pk{\cap_{j=1}^n\{\widetilde Z_{m,k}^{(\kappa)}(t_j){+x}\le z_j\}}, \quad u\to \IF
 \end{align}
 holds for all  $x > 0$ and $z_j\in\R, 1\le j\le n$.
Define below
$$
 \Delta_{iu}(t_j)= X_i(q_\kappa t_j) - r_i(q_\kappa t_j)X_i(0), \ \ \ 1\le i\le m+k,\ 1\le j\le n.
$$
By \eqref{corrr} we have
\begin{align*}
\lefteqn{u^{{2\tau/\kappa}} \Cov(\Delta_{iu}(s),\Delta_{iu}(t)) \to   C_i (s^{\alpha} +
t^{\alpha} - \abs{s-t}^\alpha )}\notag \\
&\qquad= 2C_i\Cov(Z_i(s), Z_i(t)), \quad u\to\IF, \ s,t >0, \  1\le i\le  m+k.
\end{align*}
Therefore\cL,
$$\{u^{{\tau}/\kappa}\Delta_{iu}(t), t\ge 0\}\todis\{\sqrt{2C_i}Z_i(t),t\ge 0\},\ \ u \to\IF,\ \ 1\le i\le  m+k.
$$
{Furthermore, by the independence of { $\Delta_{iu}(t)$'s }
 and $X_{i}(0)$'s, the random processes $Z_i$'s can be chosen such that they are independent of $\zeta_{m,k}^{(\kappa)}(0)$.  \ccL{Note that  $\vk X^{(1)}(0) \stackrel d= R_1 \vk{O}_1$ holds for some $R_1>0$ which is independent of $\vk{O}_1$.} Then, 
 \cL{using the Taylor\rq{}s expansion of $(1+x)^{\kappa/2} = 1+ \kappa x/2 + o(x), x\to0$}
 for any $z_j\in\R$, $1\le j\le n$ we have
\begin{align}\label{conLD: X1}
p_0(u)&= \pk{\bigcap_{j=1}^n\left\{|\vk X^{(1)}(q_\kappa t_j)|^\kappa\le {u_{\kappa, z_j}}\right\}\Big\lvert |\vk X^{(1)}(0)|^\kappa= u_{\kappa,x}}\nonumber \\
&=\mathbb P\left\{\bigcap_{j=1}^n\left\{w_\kappa(u)\left(R_1^\kappa\left(1+\frac1{R_1^2} \kE{V_u(\ccL{t_j})}\right)^{\kappa/2}-R_1^\kappa\right) \right.\right. \left.  \le z_j-x\ccL{\Bigg\}}\Bigg\lvert R_1^\kappa= u_{\kappa, x} \right\} \nonumber\\
&=\mathbb P\left\{\bigcap_{j=1}^n\left\{\frac\kappa2w_\kappa(u)R_1^{\kappa-2} \kE{V_u(t_j)}(1+o_p(1)) \le z_j-x\right\}\right.
  \Bigg\lvert R_1^\kappa= u_{\kappa, x} \ccL{\Bigg\}} \nonumber\\
 &= \pk{\bigcap_{j=1}^n\left\{\sum_{i=1}^m  \frac{\sqrt{2C_i}O_iZ_i(t_j)}{u^{(\tau-1)/\kappa}}(1+o_p(1))-
\left(\sum_{i=1}^m\frac{C_iO_i^2}{u^{2(\tau-1)/\kappa}}\right) t_j^\alpha (1+o_p(1)) +x \le
z_j\right\}},\quad u\to \IF,
\end{align}
where  $V_u(t_j):=\sum_{i=1}^m\Delta_{iu}^2(t_j) +2\sum_{i=1}^m\Delta_{iu}(t_j)r_i(q_\kappa t_j)X_i(0)-\sum_{i=1}^m(1-r^2_i(q_\kappa t_j))X_i^2(0)$.
Consequently, the claim for $k=0$ follows. Next, for $k\ge1$, we {rewrite $p_k(u)$ as} (recall that $u_{\kappa, x} = u+ x/ w_\kappa(u)$)
\begin{align}\label{conLD: XX}
\tE{p_k(u)}&= \int_0^\IF\pk{\bigcap_{j=1}^n\left\{\zeta_{m,k}^{(\kappa)}(q_\kappa t_j)\le  u_{\kappa, z_j} \right\}\Big\lvert|\vk
X^{(1)}(0)|^\kappa = u_{\kappa,x+y}, |\vk X^{(2)}(0)|^\kappa=\frac
y{w_\kappa(u)} }\nonumber\\
&\quad\times \frac{f_{|\vk X^{(1)}(0)|^\kappa}(u_{\kappa,x+y})f_{ |\vk X^{(2)}(0)|^\kappa}(y/w_\kappa(u))}
{w_\kappa(u)f_{\zeta_{m,k}^{(\kappa)}(0)}(u_{\kappa,x})}\, dy \nonumber \\
&=\int_0^\IF\pk{\bigcap_{j=1}^n\left\{|\vk X^{(1)}(q_\kappa t_j)|^\kappa\le  u_{\kappa, {z_j + w_\kappa(u)\cdot |\vk X^{(2)}(q_\kappa t_j)|^\kappa }} \right\}\Big\lvert|\vk
X^{(1)}(0)|^\kappa = u_{\kappa,x+y}, |\vk X^{(2)}(0)|^\kappa=\frac
y{w_\kappa(u)} }\times h_{\kappa, u}(y)\, dy,
\end{align}
 where 
 \BQN\label{hkappa}
 h_{\kappa, u}(y)&:=& \frac{f_{|\vk X^{(1)}(0)|^\kappa}(u_{\kappa,x+y})f_{ |\vk X^{(2)}(0)|^\kappa}(y/w_\kappa(u))}
{w_\kappa(u)f_{\zeta_{m,k}^{(\kappa)}(0)}(u_{\kappa,x})}\nonumber\\
 &=&\frac{f_{|\vk X^{(1)}(0)|^\kappa}(u_{\kappa,x+y})}{\cL{f_{|\vk X^{(1)}(0)|^\kappa}(u_{\kappa,y})}} \frac{f_{\zeta_{m,k}^{(\kappa)}(0)}(u)}{f_{\zeta_{m,k}^{(\kappa)}(0)}(u_{\kappa,x})} \frac{\cL{f_{|\vk X^{(1)}(0)|^\kappa}(u_{\kappa,y})f_{ |\vk X^{(2)}(0)|^\kappa}(y/w_\kappa(u))}}{w_\kappa(u) f_{\zeta_{m,k}^{(\kappa)}(0)}(u)} \notag\\
 &\sim& \cL{\frac{f_{ |\vk X^{(1)}(0)|^\kappa}(u_{\kappa,y}) f_{ |\vk X^{(2)}(0)|^\kappa}(y/w_\kappa(u)) }{\int_0^\IF f_{ |\vk X^{(1)}(0)|^\kappa}(u_{\kappa,y})f_{ |\vk X^{(2)}(0)|^\kappa}(y/w_\kappa(u)) \,dy}}, \quad {u\to\IF}.
 \EQN
\ccJ{Here the last step follows by \nelem{LTA}.  } 

 \ccLP{Noting that $\vk X^{(2)}(0) \stackrel d= R_2 \vk{O}_2$ holds for some $R_2>0$ which is independent of $ \vk{O}_2$,}
  we have by similar arguments as {in} \eqref{conLD: X1} that, for any $t\ge0$
\BQN \label{def_theta}
\lefteqn{ \left(w_\kappa(u)|\vk X^{(2)}(q_\kappa t)|^\kappa\right)^{2/\kappa} \Big\lvert \{w_\kappa(u)|\vk X^{(2)}(0)|^\kappa=y\}} \nonumber
\\
&&= (w_\kappa(u))^{2/\kappa} \Bigg(\sum_{i=m+1}^{m+k} \ccL{X^2_i(0)}+ 2\sum_{i=m+1}^{m+k}r_i(q_\kappa t) X_i(0) \Delta_{iu}(t)+\sum_{i=m+1}^{m+k} \Delta_{iu}^2(t)\nonumber \\
&& \quad- \sum_{i=m+1}^{m+k} (1-r_i(q_\kappa t)^2) X^2_i(0)  \Bigg)  \Bigg\lvert \left\{ R_2^\kappa= \frac y{(w_\kappa(u))^{1/\kappa}}\right\} \nonumber \\
&&=(w_\kappa(u))^{2/\kappa}\Bigg(R_2^2+ 2\frac{R_2}{u^{\tau/\kappa}}\sum_{i=m+1}^{m+k}\sqrt{2C_i}O_iZ_i(t)(1+{o_p}(1)) +\frac2{u^{2\tau/\kappa}}\sum_{i=m+1}^{m+k} C_iZ_i^2(t) (1+o_p(1))\nonumber \\
&&\quad - 2\fracl{R_2}{u^{\tau/\kappa}}^2\sum_{i=m+1}^{m+k}  C_iO_i^2 t^\alpha(1+o_p(1))\Bigg)\Bigg\lvert \left\{R_2^\kappa= \frac y{(w_\kappa(u))^{1/\kappa}}\right\}\nonumber \\
&& =y^{2/\kappa} + 2y^{1/\kappa}\fracl{w_\kappa(u)}{u^\tau}^{1/\kappa} \sum_{i=m+1}^{m+k}\sqrt{2C_i}O_iZ_i(t)(1+o_p(1))+2\fracl{w_\kappa(u)}{u^\tau}^{2/\kappa} \sum_{i=m+1}^{m+k}C_iZ_i^2(t) (1+o_p(1)) \nonumber\\
&&=: \ccJ{\kE{\theta_{\kappa,u}(y, t)}.}
\EQN
This together with \eqref{conLD: X1} and \eqref{conLD: XX} implies that
\begin{align} \label{pku}
p_k(u)&= \int_0^\IF \mathbb P\left\{\bigcap_{j=1}^n\left\{\sum_{i=1}^m  \frac{\sqrt{2C_i}O_iZ_i(t_j)}{u^{(\tau-1)/\kappa}}(1+o_p(1))-
\left(\sum_{i=1}^m\frac{C_iO_i^2}{u^{2(\tau-1)/\kappa}}\right) t_j^\alpha (1+o_p(1))+x+y \right. \right. \nonumber \\
&\quad \le
z_j +w_\kappa(u)|\vk X^{(2)}(q_\kappa t_j)|^\kappa \Bigg\}\Bigg\lvert |\vk X^{(2)}(0)|^\kappa =\frac y{w_\kappa(u)}\Bigg\}\, h_{\kappa, u}(y)\,dy \nonumber \\
&= \int_0^\IF \mathbb P\left\{ \bigcap_{j=1}^n\left\{\sum_{i=1}^m  \frac{\sqrt{2C_i}O_iZ_i(t_j)}{u^{(\tau-1)/\kappa}}(1+o_p(1))-
\left(\sum_{i=1}^m\frac{C_iO_i^2}{u^{2(\tau-1)/\kappa}}\right) t_j^\alpha (1+o_p(1))+x+y \right.\right. \nonumber \\
&\quad\le z_j  +
{(\theta_{\kappa,u}(y, t_j)})^{\frac {\kappa}{2}}\Big\} \Big\} \, h_{\kappa, u}(y)\,dy.
\end{align}
Recalling  that $\tau=\max(2/\kappa-1,1)$
 \cL{and $w_\kappa(u)=(1/\kappa) u^{2/\kappa-1}$}, we have \cL{by \eqref{def_theta} that $(\theta_{\kappa,u}(y, t_j))^{\kappa/2} =y +o_p(1)$ for \ccJ{$\kappa>1$}.  While for $\kappa\in \kE{(0,1]}$,} it follows by \ccL{\eqref{hkappa} and \nelem{LTA} that,
\BQN\label{density_g}
\cL{h_{\kappa,\infty}(y):=}\lim_{u\to\IF} h_{\kappa,u}(y)= \frac{1}{\Gamma(k/\kappa)}y^{k/\kappa-1}e^{-y},\quad \ccL{y>0},
\EQN
which is the pdf of a Gamma distributed random variable with parameter $(k/\kappa,1)$}. \cL{Hence, combining \eqref{hkappa}--
\eqref{density_g} and \eqref{KL} for the definition of $\widetilde Z_{m,k}^{(\kappa)}(t)$, the claim in  \eqref{Proof_Z} follows. Consequently, the proof of \nelem{LTA2} is complete.}  \QED

The next lemma corresponds to Condition B in  \ccJ{\cite{Albin1990}; see also
\cite{AlbinPHD, Albin2003}.}
As shown in Chapter 5  in \cite{AlbinPHD} this condition
\cL{ is crucial \ccJ{ in ensuring} that the double sum part is asymptotically \ccJ{negligible} with respect to the principal sum.} 
 \ccLP{Denote in the following by $[x]$ the integer part of $x\in\R$.}

\BEL \label{LTB} If $\{\zeta_{m,k}^{(\kappa)}(t), t\ge0\}$ is as in \netheo{T1}, then for any $T, a>0$
\BQNY \limsup_{u \to \IF} \sum_{j= N}^{[T/(aq_\kappa)]}
\pk{\zeta_{m,k}^{(\kappa)}(aq_\kappa j) > u\Big \lvert
\zeta_{m,k}^{(\kappa)}(0) > u } \to 0, \quad N\to \IF.
\EQNY
\EEL
{\bf Proof.}
{Note first that the case $k=0$ is treated in \cite{Albin1990},
p.119.} \KK{Using the fact that the standard bivariate Gaussian
distribution is exchangeable  \tE{for $u>0$}  we have
$$\pk{\zeta_{m,k}^{(\kappa)}(q_\kappa t)>u \Big\lvert  \zeta_{m,k}^{(\kappa)}(0)>u}= 2\pk{\zeta_{m,k}^{(\kappa)}(q_\kappa t)>u, |\vk X^{(1)}(q_\kappa t)| >|\vk X^{(1)}(0)|\Big\lvert  \zeta_{m,k}^{(\kappa)}(0)>u}=:\tE{ 2\Theta(u)}.$$
Further, it follows from \nelem{LTA} that, for any} $k\ge1$
\begin{align*}
\tE{\Theta(u)} &= \int_0^\IF \int_0^\IF
\pk{\zeta_{m,k}^{(\kappa)}(q_\kappa t)>u, |\vk X^{(1)}(q_\kappa
t)|>|\vk X^{(1)}(0)| \Big\lvert
|\vk X^{(1)}(0)|^\kappa= u_{\kappa,x+y}, |\vk X^{(2)}(0)|^\kappa= \frac{y}{w_\kappa(u)}} \\
& \quad \times \frac{f_{|\vk X^{(1)}(0)|^\kappa}(u_{\kappa,x+y})f_{|\vk X^{(2)}(0)|^\kappa}\fracl y{w_\kappa(u)}}{w^2_\kappa(u)\pk{\zeta_{m,k}^{(\kappa)}(0)>u}} \,dxdy \\
& \le  \int_0^\IF \int_0^\IF \pk{|\vk X^{(1)}(q_\kappa t)|^\kappa > u_{\kappa,y}\Big\lvert|\vk X^{(1)}(0)|^\kappa= u_{\kappa,x+y}}\frac{f_{|\vk X^{(1)}(0)|^\kappa}\left(u_{\kappa,x+y}\right)f_{|\vk X^{(2)}(0)|^\kappa}\fracl y{w_\kappa(u)}}{{w^2_\kappa(u)}\pk{\zeta_{m,k}^{(\kappa)}(0)>u}} \,dxdy \\
&{=}\int_0^\IF \pk{|\vk X^{(1)}(q_\kappa t)|^\kappa >
u_{\kappa,y}\Big\lvert|\vk X^{(1)}(0)|^\kappa> u_{\kappa,y}}{\frac{\pk{|\vk
X^{(1)}(0)|^\kappa>u_{\kappa,y}}}
{w_\kappa(u)\pk{\zeta_{m,k}^{(\kappa)}(0)>u}}f_{|\vk X^{(2)}(0)|^\kappa}\fracl y{w_\kappa(u)}}\,dy.
 \end{align*}
{Moreover, in
view of the treatment of the case $k=0$ in \cite{Albin1990}, p.119  we
readily see that, \ccLP{for any $p\ge1$,} with $R(t):=\max_{1\le i\le m}r_i(t), r(t):=\min_{1\le i\le m}r_i(t)$ and $\Phi(\cdot)$ denoting the $N(0,1)$ distribution function
\BQNY
\pk{|\vk X^{(1)}(q_\kappa t)|^\kappa >
u_{\kappa,y}\Big\lvert|\vk X^{(1)}(0)|^\kappa> u_{\kappa,y}}
 &\le &4m\left(1-\Phi\fracl{(1-R(q_\kappa
t)){u^{1/\kappa}}}{\sqrt{m(1-r^2{{(q_\kappa t)}})}}\right) \\
&\le & {K_p t^{-\alpha p/2},\quad \forall
q_\kappa t\in(0, T]
}
\EQNY
\ccLP{holds} for some $K_p>0$ not depending on $u, t$ and $y$. Consequently\cL,
\BQN\label{Decomp: N1}
\pk{\zeta_{m,k}^{(\kappa)}(q_\kappa t) > u \Big \lvert
\zeta_{m,k}^{(\kappa)}(0) > u} &\le& 2K_p t^{-\alpha p/2}\int_0^\IF {\frac{\pk{|\vk
X^{(1)}(0)|^\kappa>u_{\kappa,y}}}
{w_\kappa(u)\pk{\zeta_{m,k}^{(\kappa)}(0)>u}}f_{|\vk X^{(2)}(0)|^\kappa}\fracl y{w_\kappa(u)}}\,dy\notag \\
&=&2K_p t^{-\alpha p/2},\quad \forall
q_\kappa t\in(0, T].
\EQN
Therefore, with $p=4/\alpha$
\BQNY
\lefteqn{\limsup_{u \to \IF} \sum_{j= N}^{[T/(aq_\kappa)]}
\pk{\zeta_{m,k}^{(\kappa)}(aq_\kappa j) > u\Big \lvert
\zeta_{m,k}^{(\kappa)}(0) > u } }\\
&&\le 2K_p\int_{aN}^\IF x^{-2}\,dx = \frac{2K_p}{aN}\to0, \quad N\to\IF
\EQNY
establishing the proof. } \QED

\cJ{The lemma below} concerns the accuracy of the discrete approximation to the continuous process, which is related to Condition C in \cite{Albin1990}.
As shown in \cite{Albin2003} (see  Eq. (7) therein), in order to verify  Condition C  the following lemma is sufficient. \cJ{Its  proof is relegated to the \cL{appendix}.}
 \BEL \label{LTC}
If $\{\zeta_{m,k}^{(\kappa)}(t),t\ge0\}$ is  as in \netheo{T1}, then there
exist some constants $C, p>0, d>1$ and $\lambda_0, u_0 >0$ such that
\BQNY \pk{\zeta_{m,k}^{(\kappa)}(q_\kappa t) > u
+\frac{\lambda}{w_\kappa(u)}, \zeta_{m,k}^{(\kappa)}(0) \le u} \le
Ct^d\lambda^{-p} \pk{\zeta_{m,k}^{(\kappa)}(0)>u} \EQNY for $0<
t^\varpi < \lambda < \lambda_0$ and $u> u_0$. Here  \cL{${\varpi}$ is $\alpha/2$ for $\kappa\ge1$, and $(\alpha/2)\min(\kappa/(4(1-\kappa)),1)$ otherwise}. \EEL

\prooftheo{T1}
It follows from Lemmas  {\ref{LTA}--}\ref{LTC} that all the
 assumptions of Theorem 1 in
\cite{Albin1990} are satisfied by the process $\zeta_{m,k}^{(\kappa)}$, which immediately  establishes the proof. \QED

\prooftheo{T2}
{In view of \eqref{Decomp: N1} with $p=4/\alpha$} and letting $v_\kappa = v_\kappa(u) = 1/q_\kappa(u) = u^{2{\tau}/(\alpha\kappa)}${,} {we obtain }
\begin{align*}
v_\kappa\int_{N/v_\kappa}^T \pk{\zeta_{m,k}^{(\kappa)}(s) > u \Big \lvert \zeta_{m,k}^{(\kappa)}(0) > u} \, ds &  =
\int_N^{v_\kappa T} \pk{ \zeta_{m,k}^{(\kappa)}(s/v_\kappa) > u \Big \lvert \zeta_{m,k}^{(\kappa)}(0) > u}\, ds\\
&\le K_{4/\alpha}\int_N^{v_\kappa T} s^{-2}\, ds \le   \frac{K_{4/\alpha}}N,\quad u\to\IF.
\end{align*}
Hence,
\BQNY
 \lim_{N\to\IF}\limsup_{u\to\IF}v_\kappa\int_{N/v_\kappa}^T \pk{\zeta_{m,k}^{(\kappa)}(s) > u \Big \lvert \zeta_{m,k}^{(\kappa)}(0) > u} \, ds = 0.
 \EQNY
Since further \nelem{LTA2} holds, the claim follows by Theorem 3.1 in \cite{Berman82}.\QED

  As shown by Theorem 10 in \cite{Albin1990}, in order to derive
the Gumbel limit theorem for the random process
$\zeta_{m,k}^{(\kappa)}$ two {additional} conditions, which were
first addressed by the seminal \cL{contributions} 
 \cite{LeadbetterR1982,leadbetter1983extremes}, 
need to be checked, namely the mixing Condition
$D$ and the Condition $D'$ therein.
\cJi{These two
conditions will \cL{follow} from  Lemma \ref{LTD1} and \nelem{LTD2} below;}
 their proofs are displayed in the \cL{appendix}.
 \BEL \label{LTD1} Let {$T$ and $a$ be any given positive constants and $M\in(0, T)$.}
If $\{\zeta_{m,k}^{(\kappa)}(t), t\ge0\}$ is as in
 \netheo{T1},  then for any $0\le s_1 < \cdots < s_p < t_1 < \cdots < t_{p'}$ in
 $ \{ aq_\kappa j: j\in \mathbb Z, 0\le aq_\kappa j \le {T}\} $ such that $t_1 - s_p \ge M$
 \BQN \label{Asym.ind}
&& \abs{
 \pk{ \cap_{ i=1}^p\{\zeta_{m,k}^{(\kappa)}(s_i) \le u\}, \cap_{ j=1}^{p'}\{\zeta_{m,k}^{(\kappa)}(t_j) \le u\} }
 -  \pk{\cap_{ i=1}^p\{\zeta_{m,k}^{(\kappa)}(s_i) \le u\}} \pk{\cap_{ j=1}^{p'}\{\zeta_{m,k}^{(\kappa)}(t_j) \le u\}}
  }
  \notag\\
 &&\qquad
 \le K u^{\ccLP{\varsigma}} \sum_{1\le i \le p, 1\le j\le p'} {\widetilde r(t_j - s_i)} \expon{ -\frac{u^{2/\kappa}}{1+{\widetilde r(t_j - s_i)}}
 }
 \EQN
  and
 \BQN \label{Asym.indd}
&&
 \abs{
 \pk{ \cap_{ i=1}^p\{\zeta_{m,k}^{(\kappa)}(s_i) > u\}, \cap_{ j=1}^{p'}\{\zeta_{m,k}^{(\kappa)}(t_j) > u\} }
 -  \pk{\cap_{ i=1}^p\{\zeta_{m,k}^{(\kappa)}(s_i) > u\}} \pk{\cap_{ j=1}^{p'}\{\zeta_{m,k}^{(\kappa)}(t_j) > u\}}
  }\notag\\
 &&\qquad
 \le   K u^{\ccLP{\varsigma}} \sum_{1\le i \le p, 1\le j\le p'}{\widetilde r(t_j - s_i)} \expon{ -\frac{u^{2/\kappa}}{1+{\widetilde r(t_j - s_i)}}
 }
 \EQN
 hold for all $u>0$ and some $K>0$ not depending on $u$. \ccLP{Here $\varsigma =2/\kappa\cL{\big(m-k(2/\kappa-1)-1+\max(0, 2(1/\kappa-1)\big)}$} and ${\widetilde r(t):=\max_{1\le l\le m+k}|r_l(t)|}, t>0$.
 \EEL
}

 \BEL \label{LTD2}
Under the assumptions of \netheo{T3},
 for $\ccLP{\varsigma},$ {$\widetilde r(\cdot)$ as in \nelem{LTD1} and} $T_\kappa$  given by
 \BQN \label{def: T}
 T_\kappa= T_\kappa(u)= \frac{1}{{H^{m,k}_{\alpha,\kappa}} ({\vk C})} \frac{q_\kappa(u)}{\pk{\zeta_{m,k}^{(\kappa)}(0) >u}},
 \EQN
\cL{we have, }for any given constant
$\ve\in(0,T_\kappa)$
\BQN \label{Lim: 1}
 u^{\ccLP{\varsigma}} \frac {T_\kappa}{q_\kappa} \sum_{\ve\le aq_\kappa j \le T_\kappa}
{\widetilde r(aq_\kappa j)} \expon{ -\frac{u^{2/\kappa}}{1+{\widetilde r(aq_\kappa j)}} } \to
0,\quad\ u\to\IF.
 \EQN
 \EEL

\prooftheo{T3}   \cJi{To establish  Conditions $D$ and $D'$ in \cite{Albin1990}, we shall make use of \nelem{LTD1} with $T=T_\kappa$  given by \eqref{def: T} and
$M=\ve\in(0,T_\kappa)$}, and \nelem{LTD2}.  First
note that the
{right-hand} side of \eqref{Asym.ind} is bounded from above by
 $$
K u^{\ccLP{\varsigma}} \frac{T_\kappa} {aq_\kappa} \sum_{\ve\le aq_\kappa j \le T_\kappa}
{\widetilde r(aq_\kappa j)} \expon{ -\frac{u^{2/\kappa}}{1+{\widetilde r(aq_\kappa j)}}},
 $$
 which by an application of
 \eqref{Lim: 1} implies that the mixing Condition $D$ in
\cite{Albin1990} holds for the random process $\zeta_{m,k}^{(\kappa)}$.

Next, we prove Condition $D'$ in \cite{Albin1990}, i.e., for any given positive constants $a$ and ${\widetilde T}$
\BQN \label{cond: D}
 \limsup_{u\to\IF} \sum_{j =
[\widetilde T/(aq_\kappa)]}^{\left[\ve / \pk{\zeta_{m,k}^{(\kappa)}(0) >u}\right]}
\pk{\zeta_{m,k}^{(\kappa)}(aq_\kappa j)
>u \Big \lvert \zeta_{m,k}^{(\kappa)}(0)
>u} \to 0,\quad\ve\downarrow 0.
 \EQN
Indeed, \cJi{by \eqref{Asym.indd}} for some \cL{$\widetilde M>\widetilde T$} and a positive
constant $K$ \BQNY \pk{\zeta_{m,k}^{(\kappa)}(aq_\kappa j) >u \Big \lvert
\zeta_{m,k}^{(\kappa)}(0)
>u} \le \pk{\zeta_{m,k}^{(\kappa)}(0) > u} + K u^{\ccLP{\varsigma}}\frac {T_\kappa}{q_\kappa} {\widetilde r(aq_\kappa j)}
\expon{-\frac{u^{2/\kappa}}{1+{\widetilde r(aq_\kappa j)}}}
  \EQNY
   holds for 
   \cL{$u>0$}  and $aq_\kappa j>\widetilde M$. \kE{Consequently,}
\BQNY
\lefteqn{
 \limsup_{u\to\IF} \sum_{j =
[\widetilde T/(aq_\kappa)]}^{\left[\ve / \pk{\zeta_{m,k}^{(\kappa)}(0) >u}\right]}
\pk{\zeta_{m,k}^{(\kappa)}(aq_\kappa j)
>u \Big \lvert \zeta_{m,k}^{(\kappa)}(0)
>u}
}
\\
&& \le \limsup_{u\to\IF} \sum_{j = [\widetilde T/(aq_\kappa)]}^{[\widetilde M/(aq_\kappa)]}
\pk{\zeta_{m,k}^{(\kappa)}(aq_\kappa j)
>u \Big \lvert \zeta_{m,k}^{(\kappa)}(0)
>u} + \ve \notag \\
&& \quad + \limsup_{u\to\IF}  K u^{\ccLP{\varsigma}}\frac{T_\kappa}{q_\kappa}\sum_{j =
[\widetilde M/(aq_\kappa)]}^{\left[\ve / \pk{\zeta_{m,k}^{(\kappa)}(0) >u}\right]} {\widetilde r(aq_\kappa j)}
\expon{-\frac{u^{2/\kappa}}{1+{\widetilde r(aq_\kappa j)}}}, \notag
\EQNY
\cL{which equals $\epsilon$} by an
application of \nelem{LTB} and \eqref{Lim: 1},
respectively.
\ccLP{It follows then   that   \eqref{cond: D} holds}. Consequently, in view of {Theorem 10 in \cite{Albin1990}}
 we have, for $T_\kappa$ given  by \eqref{def: T} \BQNY
\limit{u} \pk{\sup_{t\in[0, T_\kappa]} \zeta_{m,k}^{(\kappa)}(t) \le u + \frac x{w_\kappa(u)}} = \expon{-e^{-x}}, \quad
x\inr.
 \EQNY
\cJ{Expressing $u$ in term{s} of $T_\kappa$ using \eqref{def: T}} {(see also \eqref{Asym: T})} we obtain
the required claim with {$a_T^{(\kappa)}, b_T^{(\kappa)}$} given by
\eqref{Normalized constant} \cJi{for any $x\in\R$; the uniform convergence in $x$ follows since all functions (with respect to $x$) are
continuous, bounded and increasing.} \QED

\section{Appendix}\label{sec4}

\prooflem{LTC}
By \eqref{corrr}, for any small  $\epsilon \in(0,1)$ there exists some
positive constant $B$  such that
\BQNY
r_i(t) \ge \frac12 \quad{\rm and}\quad
1-r_i(t) \le Bt^\alpha,\quad\forall t\in (0, \epsilon], \ 1\le i\le m+k.
 \EQNY
Furthermore, for any positive $t$ satisfying (recall $\ccL{\varpi}=\alpha/2\I{\kappa\ge1} + \alpha/2\min(\kappa/(4(1-\kappa)),1)\I{0<\kappa<1}$)
  \ccL{$$0 < t^\varpi < \lambda < \lambda_0 :=
\min\left(\frac1{2^{\kappa+4}B}, \frac{\kappa}{2^{\kappa+2}}, \epsilon^\varpi\right)$$} and {any} $u >\ccL2$
 \BQN \label{r.ineq}
 {u^{2{\tau}/ \kappa}\theta_{\kappa}(t)\le 2^\kappa \kappa B t^\alpha \le
\frac{\kappa t^{\alpha  /2}}{\ccL{16}}  \quad{\rm with}\quad \theta_\kappa(t):=\frac1{(r(q_\kappa t))^\kappa}
-1}, \ r(t) := \min_{1\le i\le {m+k}} r_i(t).
\EQN
Let
\ccL{$ (\X^{(1)}_{1/r}( t), \X^{(2)}_{1/r}(t) ):=
\big(X_1( t) - r^{-1}_1( t) X_1(0), \ldots,
X_{m+k}(t) - r^{-1}_{m+k}( t) X_{m+k}(0)\big)$
\tE{which} by definition is independent of $\{\zeta_{m,k}^{(\kappa)}(t),
t\ge0\}$}.  For  $j=1,2$
\BQN\label{D.ineq}
\pk{|\X^{(j)}_{1/r}(q_\kappa t)|>x} \le \pk{|\X^{(j)}(0)|>\frac
x{2\sqrt{2Bu^{-2{\tau}/\kappa} t^\alpha}}}, \quad
u\theta_\kappa(t)\le \frac\lambda{2w_\kappa(u)}.
\EQN
\cJi{In the following, the cases $\kappa=1, \kappa\in(1,\IF)$  and $\kappa\in(0,1)$ will be considered in turn.}

\underline{Case $\kappa=1$}: \ccLP{Note by the triangular inequality}
 \begin{align*}
\zeta_{m,k}^{(1)}(q_1 t)
&\le  |\vk X^{(1)}_{1/r}(q_1 t)|+ |\vk X^{(2)}_{1/r}(q_1
t)| + \frac1{r(q_1 t)}\zeta_{m,k}^{(1)}(0) +\theta_1(t) |\vk
X^{(2)}(0)|.
\end{align*}
Consequently, \ccLP{from  \eqref{D.ineq} we get}
\BQNY
\notag&& \pk{\zeta_{m,k}^{(1)}(q_1 t) > u +\frac{\lambda}{u},  \zeta_{m,k}^{(1)}(0) \le u} \\
\notag&& \le \pk{|\X^{(1)}_{1/r}(q_1 t)|+ |\X^{(2)}_{1/r}(q_1 t)|+ \theta_1(t)|\vk X^{(2)}(0)| > \frac{\lambda}{2u}, \zeta_{m,k}^{(1)}(q_1 t) > u} \\ 
\notag& &\le \pk{|\X^{(1)}_{1/r}(q_1 t)|+
|\X^{(2)}_{1/r}(q_1 t)| > \ccL{\frac{\lambda}{3u}}}\pk{ \zeta_{m,k}^{(1)}(q_1 t) >
u}+\pk{\theta_1(t)|\vk X^{(2)}(0)| >
\frac{\lambda}{\ccL6u}} \\
&&=: I_{1u}+I_{2u}.
\EQNY
  By \eqref{r.ineq} and
\eqref{D.ineq}, we have\cL, \tE{for any} $p>1$
 \BQNY
 \pk{|\X^{(1)}_{1/r}(q_1 t)|>
\frac{\lambda}{6u}}\le \pk{|\X^{(1)}(0)|>
\frac{\lambda}{12\sqrt{2B}t^{\alpha/2}}}\le
K\fracl{\lambda}{t^{\alpha/2}}^{-p}
\EQNY
holds with some  $K>0$
(the values of $p$ and $K$ might change from
line to line below). Similarly,
 \BQNY \pk{|\X^{(2)}_{1/r}(q_1
t)|> \frac{\lambda}{6u}}  \le K
\left(\frac{\lambda}{t^{\alpha/2}}\right)^{-p}
\EQNY
\kE{and hence }
\BQN\label{Ineq.I1} I_{1u}\le K
\left(\frac{\lambda}{t^{\alpha/2}}\right)^{-p}\pk{\zeta_{m,k}^{(1)}(0)
> u}.
\EQN
Moreover, in view of \nelem{LTA} and  \eqref{r.ineq} we have for sufficiently large $u$
that
\BQN\label{Ineq.I12}
\ \ I_{2u}  \le \frac{\pk{|\vk X^{(2)}(0)|> \frac{2\lambda u}{t^{\alpha/2}}}}{\pk{\zeta_{m,k}^{(1)}(0) > u}}\pk{\zeta_{m,k}^{(1)}(0) > u}
\le K \fracl{\lambda}{t^{\alpha/2}}^{-(p-k+2)}
u^{-(p+m-2k)}\pk{\zeta_{m,k}^{(1)}(0) > u}.
\EQN
\kE{Hence}, the claim for $\kappa=1$ follows from \eqref{Ineq.I1} and \eqref{Ineq.I12} by choosing  $p>\max(4/\alpha+k, 2k)$. \\
\underline{Case $\kappa\in(1,\IF)$}:  Denote below by
$ (\Y^{(1)}( t), \Y^{(2)}(t) ):=
\big( r^{-1}_1(t) X_1(0), \ldots,
r^{-1}_{m+k}(t) X_{m+k}(0)\big)$.
Note that $|\Y^{(1)}( t)|\le |\X^{(1)}(0)|/r(t)$ and $|\X^{(2)}(0)|\le |\Y^{(2)}(t)| \le |\X^{(2)}(0)|/r(t)$ \ccLP{for all $t<\ve$}, and for some constants $K_1, K_2>0$ whose values might change from line to line below
\BQNY\label{Ineq.triangle}
|1+x|^\kappa \ge 1+\kappa x, \quad x\in\R \quad
{\rm and}\quad (1+x)^\kappa \le 1+ K_1x+K_2 x^\kappa, \quad x\ge0.
\EQNY
We have further by the triangle inequality
\BQNY
\zeta_{m,k}^{(\kappa)}(q_\kappa t) &\le &\Big(|\Y^{(1)}(q_\kappa t)|+ |\X^{(1)}_{1/r}(q_\kappa
t)| \Big)^\kappa - \abs{ |\Y^{(2)}(q_\kappa t)| -  |\X^{(2)}_{1/r}(q_\kappa
t)|  }^\kappa \\
&\le & |\Y^{(1)}(q_\kappa t)|^\kappa  + K_1|\X^{(1)}_{1/r}(q_\kappa
t)| |\Y^{(1)}(q_\kappa t)|^{\kappa-1}+K_2|\X^{(1)}_{1/r}(q_\kappa
t)|^\kappa \\
&&-|\Y^{(2)}(q_\kappa t)|^\kappa  + \kappa |\X^{(2)}_{1/r}(q_\kappa
t)|  |\Y^{(2)}(q_\kappa t)|^{\kappa-1} \\
&\le & K_1|\X^{(1)}_{1/r}(q_\kappa
t)| |\X^{(1)}(0)|^{\kappa-1} + K_2|\X^{(1)}_{1/r}(q_\kappa
t)|^\kappa \\
&&+ K_3|\X^{(2)}_{1/r}(q_\kappa
t)| |\X^{(2)}(0)|^{\kappa-1} + \frac{\zeta_{m,k}^{(\kappa)}(0)}{(r(q_\kappa t))^\kappa} + \theta_\kappa(t) |\X^{(2)}(0)|^\kappa
\EQNY
holds for $q_\kappa t\le \epsilon$ and some constant $K_3>0$. Therefore,
with {$\mu=1/(2(\kappa-1))$ and $\varphi= \alpha/(4(\kappa-1))$},
\BQN\label{Decomp12}
\notag&& \pk{\zeta_{m,k}^{(\kappa)}(q_\kappa  t) > u +\frac{\lambda}{w_\kappa(u)},  \zeta_{m,k}^{(\kappa)}(0) \le u} \\
\notag&& \le \pk{|\X^{(1)}(0)| > \frac{\lambda^\mu u^{1/\kappa}}{t^{\varphi}}} + \pk{|\X^{(2)}(0)| > \frac{\lambda^\mu u^{1/\kappa}}{t^{\varphi}}} \\
\notag&&\quad+ \mathrm{P}\left\{ K_1|\X^{(1)}_{1/r}(q_\kappa
t)| \fracl{\lambda^\mu u^{1/\kappa}}{t^{\varphi}}^{\kappa-1} + K_2|\X^{(1)}_{1/r}(q_\kappa
t)|^\kappa + K_3|\X^{(2)}_{1/r}(q_\kappa
t)| \fracl{\lambda^\mu u^{1/\kappa}}{t^{\varphi}}^{\kappa-1}\right. \\
\notag&&\quad\left. + \theta_\kappa(t) |\X^{(2)}(0)|^\kappa \ge \frac\lambda{2w_\kappa(u)}, \zeta_{m,k}^{(\kappa)}(q_\kappa t)>u
\right\}\\
 &&=: \tilde I_{1u}+ \tilde I_{2u}+ \tilde I_{3u}.
\EQN
Note by \eqref{r.ineq}  that  ${\lambda^\mu}/{t^\varphi}>1$.
Similar arguments as in \eqref{Ineq.I12} yield that
 \BQNY
\begin{array}{l}
 \tilde I_{1u} \le K\fracl{\lambda^\mu}{t^{\varphi}}^{-(p-m+2)}u^{-(p-k(2/\kappa-1)\I{\kappa\le2})/\kappa}\pk{\zeta_{m,k}^{(\kappa)}(0) >u}\\
\tilde I_{2u} \le  K\fracl{\lambda^\mu}{t^{\varphi}}^{-(p-k+2)}u^{-(p-k+m-k(2/\kappa-1)\I{\kappa\le2})/\kappa}\pk{\zeta_{m,k}^{(\kappa)}(0) >u}
\end{array}
\EQNY
and
\BQNY
 \tilde I_{3u} &\le&\left( \pk{K_1|\X^{(1)}_{1/r}(q_\kappa
t)| \fracl{\lambda^\mu u^{1/\kappa}}{t^{\varphi}}^{\kappa-1} >\frac\lambda{8w_\kappa(u)}} + \pk{K_2|\X^{(1)}_{1/r}(q_\kappa
t)|^\kappa> \frac\lambda{8w_\kappa(u)}} \right. \\
&&\left.+ \pk{K_3|\X^{(2)}_{1/r}(q_\kappa
t)| \fracl{\lambda^\mu u^{1/\kappa}}{t^{\varphi}}^{\kappa-1}> \frac\lambda{8w_\kappa(u)}} \right) \pk{\zeta_{m,k}^{(\kappa)}(q_\kappa t)>u}\\
&&+\pk{ \theta_\kappa(t) |\X^{(2)}(0)|^\kappa > \frac\lambda{8w_\kappa(u)}}\\
&=:& (I\!I_{1u}+I\!I_{2u}+I\!I_{3u})\pk{\zeta_{m,k}^{(\kappa)}(0)>u}+ I\!I_{4u}.
\EQNY
Furthermore,
\BQNY
I\!I_{1u}&\le& \pk{|\X^{(1)}(0)|>K_1\frac{\lambda^{1/2}u^{-1/\kappa}}{t^{-\alpha/4}(r^{-2}(q_\kappa t)-1)^{1/2}} }\\
&\le& \pk{|\X^{(1)}(0)|>K_1\frac{\lambda^{1/2}}{t^{\alpha/4}} }
\le K\fracl{\lambda}{t^{\alpha/2}}^{-p/2}.
\EQNY
Similarly,
\BQNY
I\!I_{2u}\le K\fracl{\lambda u^{2(1-1/\kappa)}}{t^{\alpha\kappa/2}}^{-p/\kappa}, \quad I\!I_{3u}\le K\fracl{\lambda}{t^{\alpha/2}}^{-p/2}.
\EQNY
Next, we deal with $I\!I_ {4u}$. We have by \eqref{r.ineq} that $2^{\kappa+4}B t^{\alpha/2}\le 1$.
\mE{Hence as in the proof of}  \eqref{Ineq.I12}\cL{, we have}
\BQN\label{I4u}
I\!I_{4u}&\le& \pk{|\X^{(2)}(0)|^\kappa > \frac{2\lambda u}{t^{\alpha/2}}\frac1{2^{\kappa+4}B t^{\alpha/2}}}\notag \\
&\le& K\fracl{\lambda}{t^{\alpha/2}}^{-(p-k+2)}u^{-(p-k+m-k(2/\kappa-1)\I{\kappa\le2})/\kappa}\pk{\zeta_{m,k}^{(\kappa)}(0) >u}.
\EQN
Therefore, the claim for $\kappa\in(1,\IF)$ follows from \eqref{Decomp12} and the inequalities for $\tilde I_{1u}, \tilde I_{2u}$ \cL{and $I\!I_{1u}$-- $I\!I_{4u}$}
by choosing  { $p>\max(8(\kappa-1)/\alpha+k+m, 2k)$}.

\underline{Case $\kappa\in(0,1)$}:  Note that
\BQNY
\cL{(1+x)^\kappa \le 1+x, \quad x\ge0 \quad
{\rm and}\quad -|1-x|^\kappa \le -(1-x), \quad {x\in[0,\IF)}.}
\EQNY
We have further by the triangle inequality
\BQNY
\zeta_{m,k}^{(\kappa)}(q_\kappa t) &\le &\Big(|\Y^{(1)}(q_\kappa t)|+ |\X^{(1)}_{1/r}(q_\kappa
t)| \Big)^\kappa - \abs{ |\Y^{(2)}(q_\kappa t)| -  |\X^{(2)}_{1/r}(q_\kappa
t)|  }^\kappa \\
&\le & |\Y^{(1)}(q_\kappa t)|^\kappa  + |\X^{(1)}_{1/r}(q_\kappa
t)| |\X^{(1)}(0)|^{\kappa-1}-|\Y^{(2)}(q_\kappa t)|^\kappa  + |\X^{(2)}_{1/r}(q_\kappa
t)|  |\Y^{(2)}(q_\kappa t)|^{\kappa-1} \\
&\le & \frac{|\X^{(1)}(0)|^\kappa}{(r(q_\kappa t))^\kappa}  +|\X^{(1)}_{1/r}(q_\kappa
t)| |\X^{(1)}(0)|^{\kappa-1}-|\X^{(2)}(0)|^\kappa  + |\X^{(2)}_{1/r}(q_\kappa
t)|  |\X^{(2)}(0)|^{\kappa-1} \\
&= & \frac{\zeta_{m,k}^{(\kappa)}(0)}{(r(q_\kappa t))^\kappa} + \theta_\kappa(t) |\X^{(2)}(0)|^\kappa +  |\X^{(1)}_{1/r}(q_\kappa
t)| |\X^{(1)}(0)|^{\kappa-1} + |\X^{(2)}_{1/r}(q_\kappa
t)| |\X^{(2)}(0)|^{\kappa-1}.
\EQNY
Therefore, we have by \eqref{D.ineq}, \ccL{with $\psi=\alpha/(4(1-\kappa))$}
\BQNY
\lefteqn{\pk{\zeta_{m,k}^{(\kappa)}(q_\kappa t) > u+\frac{\lambda}{w_\kappa(u)}, \zeta_{m,k}^{(\kappa)}(0) \le u}} \\
&&\le \pk{\theta_\kappa(t) |\X^{(2)}(0)|^\kappa + \frac{|\X^{(1)}_{1/r}(q_\kappa
t)|}{(u^{-\tau/\kappa}t^\psi)^{1-\kappa}} + \frac{|\X^{(2)}_{1/r}(q_\kappa
t)|}{(u^{-\tau/\kappa}t^\psi)^{1-\kappa}}> \frac{\lambda}{2w_\kappa(u)}, \zeta_{m,k}^{(\kappa)}(q_\kappa t) > u} \\
&&\quad+\pk{|\X^{(1)}(0)| \le u^{-\frac\tau\kappa}t^\psi, \zeta_{m,k}^{(\kappa)}(q_\kappa t) > u} + \pk{|\X^{(2)}(0)| \le u^{- \frac\tau\kappa}t^\psi, \zeta_{m,k}^{(\kappa)}(0) \le u, \zeta_{m,k}^{(\kappa)}(q_\kappa t) > u+\frac{\lambda}{w_\kappa(u)}}\\
&&=: I^*_{1u}+ I^*_{2u} +  I^*_{3u}.
\EQNY
Now we deal with the three terms one by one. Clearly, for any $u>2$
\BQNY
 I^*_{1u}&\le& \pk{\theta_\kappa(t) |\X^{(2)}(0)|^\kappa >\frac{\lambda}{6w_\kappa(u)}} + \pk{|\X^{(1)}_{1/r}(q_\kappa
t)| > \frac{\lambda \kappa t^{\alpha/4}}{6u^{\tau/\kappa}}}\pk{\zeta_{m,k}^{(\kappa)}(q_\kappa t) > u} \\
&&+\pk{|\X^{(2)}_{1/r}(q_\kappa
t)| > \frac{\lambda \kappa t^{\alpha/4}}{6u^{\tau/\kappa}}}\pk{\zeta_{m,k}^{(\kappa)}(q_\kappa t) > u},
\EQNY
where the first term \ccLP{can be treated as  {for} $I\!I_{4u}$, see \eqref{I4u}.} For the rest two terms, {we have, by } using \eqref{D.ineq}
\BQN\label{I1*}
\pk{|\X^{(j)}_{1/r}(q_\kappa
t)| > \frac{\lambda \kappa t^{\alpha/4}}{6u^{\tau/\kappa}}}\le \pk{|\X^{(j)}(0)| >\frac\kappa{12\sqrt{2B}}\frac{\lambda}{t^{\alpha/4}}} \le K\fracl{\lambda}{t^{\alpha/4}}^{-p}, \quad j=1,2.
\EQN
In order to deal with $I^*_{2u}$ and $I^*_{3u}$, set below
\ccL{$ (\X^{(1)}_{r}(t), \X^{(2)}_{r}(t) ):=
\big(X_1(0) - r_1( t) X_1( t), \ldots,
X_{m+k}(0) - r_{m+k}( t) X_{m+k}(t)\big)$
\tE{which} by definition is independent of $\{\zeta_{m,k}^{(\kappa)}(t),
t\ge0\}$}.  For  $j=1,2$
\BQN\label{D.ineq1}
\pk{|\X^{(j)}_{r}(q_\kappa t)|>x} \le \pk{|\X^{(j)}(0)|>\frac
{2 \sqrt \lambda x}{\sqrt{u^{-2{\tau}/\kappa} t^\alpha}}}.
\EQN
Using \ccL{further} the triangle inequality $ |\X^{(1)}_{r}(q_\kappa t)|^\kappa\ge (r(q_\kappa t))^\kappa |\X^{(1)}(q_\kappa t)|^\kappa-|\X^{(1)}(0)|^\kappa $ \ccL{and \eqref{r.ineq} (recalling} $|\X^{(1)}(q_\kappa t)|^\kappa\ge \zeta_{m,k}^{(\kappa)}(q_\kappa t)>u $), \cL{we have}
\BQN\label{I2*}
I^*_{2u} &\le& \pk{ |\X^{(1)}_{r}(q_\kappa t)|^\kappa >u\left((r(q_\kappa t))^\kappa - \frac{t^{\psi\kappa}}{u^{1+\tau}}\right)} \pk{\zeta_{m,k}^{(\kappa)}(q_\kappa t) > u}
\notag\\
&\le& \pk{ |\X^{(1)}_{r}(q_\kappa t)|^\kappa  >\frac{(1-2^{-\kappa})u}{2^\kappa}} \pk{\zeta_{m,k}^{(\kappa)}(0) > u}
\notag\\
&\le& \pk{ |\X^{(1)}(0)| >(1-2^{-\kappa})^{1/\kappa}\frac{\sqrt\lambda}{t^{\alpha/2}}} \pk{\zeta_{m,k}^{(\kappa)}(0) > u}\notag\\
&\le & K \fracl{\lambda}{t^{\alpha}}^{-p/2} \pk{\zeta_{m,k}^{(\kappa)}(0) > u}.
\EQN
For $I^*_{3u}$, {using} 
  $|\X^{(1)}(q_\kappa t)|^\kappa> u +\lambda/w_\kappa(u)$ and
\BQNY
|\X^{(1)}(0)|^\kappa =  \zeta_{m,k}^{(\kappa)}(0) + |\X^{(2)}(0)|^\kappa  \le u \left(1+\frac{t^{\psi\kappa}}{u^{1+\tau}}\right)
\EQNY
we have
\BQNY
I^*_{3u} &\le& \pk{ |\X^{(1)}_{r}(q_\kappa t)|^\kappa >u\left((r(q_\kappa t))^\kappa \left(1+\frac\lambda{uw_\kappa(u)}\right)-\left(1+\frac{t^{\psi\kappa}}{u^{1+\tau}}\right)\right)} \pk{\zeta_{m,k}^{(\kappa)}(q_\kappa t) > u} \\
&=& \pk{ |\X^{(1)}_{r}(q_\kappa t)|^\kappa >u^{-\tau}\left(\lambda \kappa(r(q_\kappa t))^\kappa - u^{1+\tau}(1-(r(q_\kappa t))^\kappa) - t^{\psi\kappa}\right)} \pk{\zeta_{m,k}^{(\kappa)}(0) > u},
\EQNY
where by \eqref{r.ineq}
\BQNY
\lambda \kappa(r(q_\kappa t))^\kappa - u^{1+\tau}(1-(r(q_\kappa t))^\kappa) - t^{\psi\kappa}
 \ge \frac{\lambda\kappa}{2^{\kappa+1}} - t^{\psi\kappa} \ge \frac{\lambda\kappa}{2^{\kappa+2}}.
\EQNY
Consequently, \cL{it follows further by \eqref{D.ineq1}} that
\BQNY
I^*_{3u} &\le& \pk{ |\X^{(1)}(0)|>2^{-2/\kappa}\kappa^{1/\kappa} \frac{\lambda^{1/\kappa+1/2}}{t^{\alpha/2}}}\pk{\zeta_{m,k}^{(\kappa)}(0) > u} \\
&\le& K\fracl{\lambda^{1/\kappa+1/2}}{t^{\alpha/2}}^{-p}\pk{\zeta_{m,k}^{(\kappa)}(0) > u},
\EQNY
which together with \eqref{I4u}, \eqref{I1*} and \eqref{I2*} completes the proof for $\kappa\in(0,1)$ by taking $p>4/\alpha+k$.
Consequently, the desired claim of \nelem{LTC} follows. This completes the  proof. \QED

\prooflem{LTD1} \cJi{We give only  the proof for \eqref{Asym.ind} since  \eqref{Asym.indd}  follows  by similar arguments.}
Since the claims for $k=0$ are already shown  in \cite{AlbinPHD}, we only consider that $k\ge1$ below.
  Define, for $j=1,2,$
independent random vectors $\left( |\Y^{(j)}(s_1)|, \ldots,
|\Y^{(j)}(s_p)| \right) $ and
$\Big(|\widetilde{\Y}^{(j)}(t_1)|, \ldots,
|\widetilde\Y^{(j)}(t_{p'})| \Big)$, {which are independent of the process $\zeta_{m,k}^{(\kappa)}$ and
have the same} distributions as those of $\left( |\X^{(j)}(s_1)|, \ldots,
|{\X^{(j)}(s_p)}| \right)$ and $\left( |{\X^{(j)}(t_{1})}|,
\ldots, |{\X^{(j)}(t_{p'})}|\right)$, respectively.  Note that, for any $u>0$, the left-hand side of \eqref{Asym.ind} is clearly bounded from above by
\BQN
&&\abs{
 \pk{ \bigcap_{ i=1}^p\left\{|{\X^{(2)}(s_i)} |^\kappa\ge |\X^{(1)}(s_i)|^\kappa-u \right\}, \bigcap_{ j=1}^{p'}\left\{|{\X^{(2)}(t_j)} |^\kappa\ge |\X^{(1)}(t_j)|^\kappa-u} \right\} \right.
 \notag \\
 && \quad - \left.
 \pk{ \bigcap_{ i=1}^p\left\{|\Y^{(2)}(s_i)|^\kappa \ge |\X^{(1)}(s_i)|^\kappa-u \right\}, \bigcap_{ j=1}^{p'}\left\{|\widetilde\Y^{(2)}(t_j)|^\kappa \ge |\X^{(1)}(t_j)|^\kappa-u\right \} }
 }
 \notag \\
 && + \abs{
 \pk{ \bigcap_{ i=1}^p\left\{|{\X^{(1)}(s_i)}|^\kappa\le |\Y^{(2)}(s_i)|^\kappa+ u \right\}, \bigcap_{ j=1}^{p'}\left\{|{\X^{(1)}(t_j)} |^\kappa\le |\widetilde\Y^{(2)}(t_j)|^\kappa +u \right\} } \right.
 \notag \\
 && \quad - \left.
\pk{\bigcap_{ i=1}^p\left\{|\Y^{(1)}(s_i)|^\kappa \le |\Y^{(2)}(s_i)|^\kappa+ u \right\}, \bigcap_{ j=1}^{p'}\left\{|\widetilde\Y^{(1)}(t_j)|^\kappa \le |\widetilde\Y^{(2)}(t_j)|^\kappa +u \right\} }
 }. \label{eq: Term}
\EQN
 Next, \ccL{note \ccL{by Cauchy-Schwarz inequality} that $u^2+v^2 \le (u^2-2\rho uv+v^2)/(1-|\rho|)$ for all $\rho\in(-1,1)$ and $u,v\in\R$. It follows that{,}  $f_{ij}(\cdot,\cdot)${,} the joint density function of
$\left(|\X^{(1)}(s_i)|, |\X^{(1)}(t_j)|\right)$ {, satisfies} 
 }
{
\BQNY
f_{i,j}(x, y) &=&\int_{|\vk x|=x, |\vk y|=y}\prod_{l=1}^m \frac1{2\pi\sqrt{1-r_l^2(t_j-s_i)}} \expon{-\frac{x_l^2-2r_l(t_j-s_i)x_ly_l+y_l^2}{2(1-r_l^2(t_j-s_i))}}\,d\vk x d \vk y\\
&\le& \frac1{(2\pi)^m (1-(\widetilde r(t_j-s_i))^2)^{m/2}} \int_{|\vk x|=x, |\vk y|=y}\prod_{l=1}^m \expon{-\frac{x_l^2+y_l^2}{2(1+|r_l(t_j-s_i)|)}}\,d\vk x d \vk y\\
&\le& \frac1{(2\pi)^m (1-(\widetilde r(t_j-s_i))^2)^{m/2}} \expon{-\frac{x^2+y^2}{2(1+\widetilde r(t_j-s_i))}}\int_{|\vk x|=x, |\vk y|=y} \,d\vk x d \vk y\\
&=& \frac{(xy)^{m-1}}{2^{m-2}(\Gamma(m/2))^2 (1-(\widetilde r(t_j-s_i))^2)^{m/2}}\expon{-\frac{x^2+y^2}{2(1+\widetilde r(t_j-s_i))}}, \quad \ccL{x, y>0.}
\EQNY
}
Therefore, in view of Lemma 2 in \cite{AlbinPHD}, with
$K$  a constant \cJi{whose} value might change from line to line,
the first \cJi{absolute value} in \eqref{eq: Term} is bounded from above by
\begin{align*}
\lefteqn{K\sum_{i=1}^p \sum_{j=1}^{p'} \int_{x^\kappa>u}\int_{y^\kappa>u}{\widetilde r(t_j-s_i)} \Big((x^\kappa-u)(y^\kappa-u)\Big)^{(k-1)/\kappa}  \expon{-\frac{(x^\kappa-u)^{2/\kappa} +
(y^\kappa-u)^{2/\kappa}}{2(1+{\widetilde r(t_j-s_i)})}} f_{ij}(x, y) \,
dx dy  }\\
 & \le  K\sum_{i=1}^p \sum_{j=1}^{p'} {\widetilde r(t_j-s_i)} \Bigg( \int_u^\IF (x-u)^{(k-1)/\kappa}x^{m/\kappa -1}
 \expon{-\frac{x^{2/\kappa}}{2(1+{\widetilde r(t_j-s_i)})}}\,
dx\Bigg) ^2\qquad \qquad\qquad \qquad \\
 &\le K u^{(2/\kappa)(m-\cL{(k-1)(2/\kappa-1)}-2)}\sum_{i=1}^p \sum_{j=1}^{p'}   {\widetilde r(t_j-s_i)}\expon{ -\frac{u^{2/\kappa}}{1+{\widetilde r(t_j-s_i)}}}{,}
\end{align*}
where \cL{in the first inequality, we use first the bound $e^{-x}\le 1, x\ge0$ and then \ccJ{a change of variable} $x\rq{}=x^\kappa$,} 
while  the second inequality follows by a change of variable
$x'=\cL{u^{2/\kappa-1}}(x-u)$ and \cL{Taylor\rq{}s expansion of
$(u+x\rq{}/u^{2/\kappa-1})^{2/\kappa}= u^{2/\kappa} +(2/\kappa)x\rq{}+O(u^{-2/\kappa})$ for large $u$ and $x\rq{}\ge0$}. 
Similarly,  denoting by \pE{$g(\cdot)$}  the pdf of  $|\vk{X}^{(2)}(0)|$,  we obtain that the second \cJi{absolute value in}  \eqref{eq: Term}
is bounded  from above  by
\BQNY \lefteqn{ K\sum_{i=1}^p \sum_{j=1}^{p'}
{\widetilde r(t_j-s_i)} \int_0^\IF\int_0^\IF\Big((x^\kappa+u)(y^\kappa+u)\Big)^{(m-1)/\kappa}\expon{-\frac{(x^\kappa+u)^{2/\kappa} +
(y^\kappa+u)^{2/\kappa}}{2(1+{\widetilde r(t_j-s_i)})}} g(x) g(y) \,
dx dy}\\
&& \le K\sum_{i=1}^p \sum_{j=1}^{p'} {\widetilde r(t_j-s_i)} \Bigg(\int_0^\IF(x^\kappa+u)^{(m-1)/\kappa} x^{k-1} \expon{-\frac{
(x^\kappa+u)^{2/\kappa}}{2(1+{\widetilde r(t_j-s_i)})}} \, d x  \Bigg)^2 \qquad\qquad \qquad\qquad\\
&& \le K u^{(2/\kappa)(m-k\cL{(2/\kappa-1)}-1)} \sum_{i=1}^p \sum_{j=1}^{p'} {\widetilde r(t_j-s_i)}
\expon{-\frac{u^{2/\kappa}}{1+{\widetilde r(t_j-s_i)}}},
\EQNY
 where the last step follows by a change of variable $\cL{x'=u^{2/\kappa-1}x^\kappa}$.
Hence the proof of \eqref{Asym.ind} is established since
$$
\big(m-k(2/\kappa-1)-1\big)-\big(m-(k-1)(2/\kappa-1)-2\big)=-2(1/\kappa-1).
$$
The desired result in \nelem{LTD1} follows. \QED

 {\prooflem{LTD2}} The proof follows by the same arguments as for Lemma 12.3.1
in \cite{leadbetter1983extremes}, using alternatively the following asymptotic relation  (recall \eqref{def: T} and \nelem{LTA})
 \BQN \label{Asym: T}
 u^{2/\kappa} = 2 \pE{\ln } T_\kappa + K_0 \ln\ln T_\kappa+
 \ln D_0(1+o(1)),\quad T_\kappa\to\IF
 \EQN
with $D_0, K_0$ defined in \netheo{T3}.
We split the sum in \eqref{Lim: 1} at
$T_\kappa^\beta$, where $\beta$ is a constant such that $0< \beta <
(1-\delta)/(1+\delta)$ and $\delta = \sup\{\widetilde r(t): t\ge
\epsilon\} <1$ (\cL{see, e.g.,} Lemma 8.1.1 (i) in
\cite{leadbetter1983extremes}). Below $K$ is again a positive
constant which {value} might change from line to line. From \eqref{Asym: T}
we conclude that  $\expon{-u^{2/\kappa}/2} \le K/T_\kappa$ and $u^{2/\kappa} = 2 \pE{\ln }T_\kappa
(1+o(1))$. Further,
\BQNY \lefteqn{ u^{\ccLP{\varsigma}} \frac{T_\kappa}{q_\kappa} \sum_{\ve\le aq_\kappa j \le T_\kappa^\beta}\widetilde r(aq_\kappa j)\expon{
-\frac{u^{2/\kappa}}{1+\widetilde r(aq_\kappa j)}}
}\\
&& \le u^{ \ccLP{\varsigma}+  \frac{4\tau}{\alpha \kappa}}  T_\kappa^{\beta+1}
\expon{-\frac{u^{2/\kappa}}{1+\delta}}
 \le K ( \pE{\ln } T_\kappa)^{ \ccL{\frac{\kappa \ccLP{\varsigma}}2} + \frac{2\tau}\alpha} T_\kappa^{\beta+1-\frac{2}{1+\delta}},
 \EQNY
 which tends to 0 as  $T_\kappa\to\IF$ since $\beta+1 - 2/(1+\delta) <0$. For the remaining sum, denoting $\delta(t)=\sup\{\abs{\widetilde r(s) \pE{\ln } s}: s\ge
 t\}$, $t>0$, we have $\widetilde r(t) \le \delta(t) / \pE{\ln } t$ as $t\to\IF$, and
 thus in view of \eqref{Asym: T} for $aq_\kappa j \ge T_\kappa^\beta$
 \BQNY
 \lefteqn{ \expon{-\frac{u^{2/\kappa}}{1+\widetilde r(aq_\kappa j)}} \le \expon{-u^{2/\kappa}\left(1-\frac{\delta(T_\kappa^\beta)}{ \pE{\ln }  T_\kappa^\beta}\right)} } \\
  && \le K \exp(-u^{2/\kappa})  \le K T_\kappa^{-2} ( \pE{\ln } T_\kappa)^{-K_0}.
 \EQNY
 Consequently, \cL{with $c$ given by \netheo{T3}} \cLL{(recall $\tau=2\max(1/\kappa-1, 0)+1$)},
 \BQN
\lefteqn{ u^{\ccLP{\varsigma}} \frac{T_\kappa}{q_\kappa} \sum_{T_\kappa^\beta \le aq_\kappa j \le T_\kappa} \ccL
{\widetilde r(aq_\kappa j)} \expon{ -\frac{u^{2/\kappa}}{1+\widetilde r(aq_\kappa j)} } } \notag \\
&& \le K u^{\ccLP{\varsigma}} \fracl{T_\kappa}{q_\kappa}^2
T_\kappa^{-2}( \pE{\ln } T_\kappa)^{-K_0}\frac{1}{ \cL{\big(\pE{\ln } T_\kappa^\beta\big)^c}}
\frac{1}{{T_\kappa}/{q_\kappa}}\sum_{T_\kappa^\beta \le aq_\kappa j \le T_\kappa} \widetilde r(aq_\kappa j) \cL{\big(\pE{\ln } (aq_\kappa j) \big)^c}
\notag \\
&& \le K (\ln T_\kappa)^{\ccL{\frac{\kappa \ccLP{\varsigma}}2} +\frac{2\tau}\alpha-K_0-\cL{c}} \frac{1}{{T_\kappa}/{q_\kappa}}\sum_{T_\kappa^\beta \le aq_\kappa j \le T_\kappa} \widetilde r(aq_\kappa j)\cL{\big(\pE{\ln } (aq_\kappa j) \big)^c}.\label{Lim: 2} 
\EQN
Since $K_0=m-2+2\tau/\alpha+k\min(1-2/\kappa,0)$ and $\kE{\ccLP{\varsigma}:=2/\kappa(m-k\cL{(2/\kappa-1)}-1+\max(0, 2\cL{(1/\kappa-1)}))}$, we have
\cL{$\ccL{{\kappa \ccLP{\varsigma}}/2} +{2\tau}/\alpha-K_0-\cL{c} = 0$ for all $\kappa>0$}. Noting further that
the Berman-type condition $\limit{t} \widetilde r(t) \cL{\big(\pE{\ln } t \big)^c}=0$ holds and $\beta<1,$ the right-hand side of
\eqref{Lim: 2} tends to 0 as
$u\to\IF$. Thus the proof is complete.
\QED

{\bf Acknowledgments}.  We would like to thank the associate editor and the two anonymous referees  for their \cccL{constructive} suggestions and corrections which greatly improved the paper.

\newcommand{\nosort}[1]{}\def\polhk#1{\setbox0=\hbox{#1}{\ooalign{\hidewidth
  \lower1.5ex\hbox{`}\hidewidth\crcr\unhbox0}}}
  \def\polhk#1{\setbox0=\hbox{#1}{\ooalign{\hidewidth
  \lower1.5ex\hbox{`}\hidewidth\crcr\unhbox0}}} \def\cprime{$'$}
  \def\cprime{$'$} \def\cprime{$'$}

\end{document}